\documentclass[a4paper]{article}
\usepackage{lmodern}
\pdfoutput=1
\usepackage{graphicx,wrapfig}
\usepackage{lipsum}
\usepackage{amsmath,amssymb,mathrsfs}
\usepackage{hyperref}
\usepackage{url}
\usepackage{mathtools}
\usepackage{changepage}
\usepackage[ruled,vlined]{algorithm2e}
\usepackage[T1]{fontenc}
\usepackage{lmodern}
\usepackage[utf8]{inputenc}
\usepackage{amsmath,amssymb,amsthm,mathrsfs}
\usepackage{graphicx,booktabs,array,longtable,enumitem,tikz}
\usepackage[section]{placeins}
\usepackage{xurl}
\usetikzlibrary{calc,positioning}
\usepackage{lmodern}
\usepackage{amsmath}
\usepackage{multirow}
\usepackage{wrapfig}
\usepackage{subfig}

\usepackage{color}
\usepackage{epsfig,amsmath,amsfonts,amssymb,amsthm,mathrsfs,ifpdf}
\usepackage{indentfirst,relsize}
\usepackage{setspace,graphicx}
\usepackage{latexsym}
\usepackage[all]{xy}
\usepackage[dvipsnames]{pstricks}
\usepackage{pst-grad} 
\usepackage{pst-plot} 
\usepackage[margin = 2.50cm]{geometry}

\usepackage{authblk} 
\usepackage{algpseudocode}
\usepackage{color,soul}
\usepackage{enumitem}

\usepackage{orcidlink}
\usepackage{tabularx}
\usepackage{booktabs,tabularx,array}
\usepackage{tikz}
\usetikzlibrary{positioning}

\usepackage[ruled,vlined]{algorithm2e}
\usepackage{booktabs,tabularx,array}
\usepackage{amssymb}
\usepackage{orcidlink}

\usepackage{tikz}
\usetikzlibrary{arrows.meta,positioning,fit,calc,backgrounds}
\usepackage{enumitem}

\newcommand{\ceil}[1]{\left\lceil #1 \right\rceil}

\newcommand{\remove}[1]{}

\newtheorem{theorem}{Theorem}
\newtheorem{lemma}[theorem]{Lemma}
\newtheorem{proposition}[theorem]{Proposition}
\newtheorem{corollary}[theorem]{Corollary}

\newtheorem{definition}[theorem]{Definition}
\newtheorem{remark}{Remark}

\usepackage{todonotes}
\usepackage[most]{tcolorbox}

\usetikzlibrary{calc,positioning}
\usepackage{orcidlink}
\hypersetup{colorlinks=true,linkcolor=blue!55!black,citecolor=blue!55!black,urlcolor=blue!55!black}
\theoremstyle{definition}
\newcommand{\rhooe}{\rho_e^o}
\newcommand{\nuI}{\nu_I}

\setlist{itemsep=3pt,topsep=5pt}
\tikzset{eop graph/.style={x=1cm,y=1cm,line width=.45pt,font=\small},
	eop vertex/.style={circle,fill=black,inner sep=1.5pt,outer sep=0pt},
	eop label/.style={font=\small,inner sep=1pt},
	eop witness/.style={line width=1.25pt}}
\newcommand{\eopv}[3]{\node[eop vertex] (#1) at (#2,#3) {};}
\newcommand{\eope}[2]{\draw (#1)--(#2);}
\newcommand{\eopfan}[3]{%
	\foreach \a in {-23,23}{\draw (#1)--++({#2+\a}:.55) node[eop vertex] {};}
	\node[eop label] at ($(#1)+({#2}:.50)$) {$\cdots$};
	\pgfmathsetmacro{\eopanchor}{#2+180}%
	\node[eop label,anchor=\eopanchor] at ($(#1)+({#2}:.68)$) {$#3$};}
\newcommand{\eoppanel}[3]{\begin{minipage}[b]{#1}\centering
		\begin{tikzpicture}[eop graph]#3\end{tikzpicture}\par\smallskip
		{\small #2}\end{minipage}}
\newcommand{\eoppath}[2]{\foreach \i in {1,...,#1}{\eopv{v\i}{.48*\i}{0}}
	\foreach \i [evaluate=\i as \j using int(\i+1)] in {1,...,#2}{\eope{v\i}{v\j}}}
\newcommand{\eoptriangle}{\eopv{u}{0}{.45}\eopv{v}{-.52}{-.4}\eopv{w}{.52}{-.4}\draw(u)--(v)--(w)--(u);}
\newcommand{\eopsquare}{\eopv{u}{0}{.55}\eopv{v}{-.55}{0}\eopv{w}{0}{-.55}\eopv{z}{.55}{0}\draw(u)--(v)--(w)--(z)--(u);}
\newcommand{\eoptheta}[1]{\eopsquare\ifnum#1=1\draw(u)--(w);\else\eopv{a}{0}{0}\draw(u)--(a)--(w);\fi}
\newcommand{\eoptriad}[3]{%
	\eopv{c}{0}{0}
	\foreach \angle/\len/\id/\cnt in {90/#1/a/r_1,210/#2/b/r_2,330/#3/d/r_3}{%
		\node[eop vertex](\id s) at ($(c)+(\angle:{.47*(\len-1)})$) {};
		\node[eop vertex](\id e) at ($(c)+(\angle:{.47*\len})$) {};
		\ifnum\len=3\node[eop vertex](\id m) at ($(c)+(\angle:.47)$) {};\draw(c)--(\id m)--(\id s);\else\draw(c)--(\id s);\fi
		\draw(\id s)--(\id e);\eopfan{\id s}{\angle+55}{\cnt}}}

\newcommand{\EOPFigureA}{%
	\begin{figure}[tbp]\centering
		\eoppanel{.31\linewidth}{$\mathcal A_1$}{\eoppath{5}{4}\eopfan{v2}{90}{a}\eopfan{v4}{90}{b}}
		\hfill\eoppanel{.31\linewidth}{$\mathcal A_2$}{\eoppath{6}{5}\eopfan{v2}{90}{a}\eopfan{v5}{90}{b}}
		\hfill\eoppanel{.31\linewidth}{$\mathcal A_3$}{\eoppath{7}{6}\eopfan{v2}{90}{a}\eopfan{v6}{90}{b}}
		\par\medskip
		\eoppanel{.23\linewidth}{$\mathcal A_4$}{\eoptriangle\eopfan{u}{90}{a}}
		\hfill\eoppanel{.23\linewidth}{$\mathcal A_5$}{\eopsquare\eopfan{u}{90}{a}}
		\hfill\eoppanel{.23\linewidth}{$\mathcal A_6$}{\eopv{c}{0}{.6}\eopv{x}{-.55}{0}\eopv{y}{.55}{0}\eopv{xx}{-.55}{-.55}\eopv{yy}{.55}{-.55}\draw(c)--(x)--(xx);\draw(c)--(y)--(yy);\eopfan{c}{90}{s-2}}
		\hfill\eoppanel{.23\linewidth}{$\mathcal A_7$}{\eopv{c}{0}{.6}\eopv{x}{0}{0}\eopv{y}{-.4}{-.55}\eopv{z}{.4}{-.55}\draw(c)--(x)--(y);\draw(x)--(z);\eopfan{c}{90}{s-1}}
		\caption{The families $\mathcal A_1,\ldots,\mathcal A_7$. A fan represents the labelled number of new leaves and is absent when that number is zero. The parameters satisfy the bounds in the definitions.}
		\label{fig:Families_graphs_A_class}
\end{figure}}

\newcommand{\EOPFigureRone}{%
	\begin{figure}[tbp]\centering
		\eoppanel{.31\linewidth}{$\mathcal R_1$(a): $K_3$, $s\ge2$}{\eopv{c}{0}{0}\eopv{x}{.7}{0}\eopv{y}{1.2}{.45}\eopv{z}{1.2}{-.45}\draw(c)--(x)--(y)--(z)--(x);\eopfan{c}{180}{s-1}}
		\hfill\eoppanel{.31\linewidth}{$\mathcal R_1$(a): $P_4$, $s\ge3$}{\eopv{c}{0}{0}\eopv{x}{.6}{0}\eopv{y}{1.1}{.25}\eopv{z}{1.6}{.25}\eopv{w}{1.05}{-.45}\draw(c)--(x)--(y)--(z);\draw(x)--(w);\eopfan{c}{180}{s-1}}
		\hfill\eoppanel{.31\linewidth}{$\mathcal R_1$(a): $K_{1,3}$, $s\ge4$}{\eopv{c}{0}{0}\eopv{x}{.7}{0}\eopv{y}{1.2}{.5}\eopv{z}{1.3}{0}\eopv{w}{1.2}{-.5}\draw(c)--(x)--(y);\draw(x)--(z);\draw(x)--(w);\eopfan{c}{180}{s-1}}
		\par\medskip
		\eoppanel{.31\linewidth}{$\mathcal R_1$(b): path at $x$\\$s\ge3$}{\eopv{c}{0}{0}\eopv{x}{.55}{.35}\eopv{y}{1.05}{.55}\eopv{z}{1.55}{.55}\eopv{w}{.55}{-.4}\eopv{a}{1.1}{-.6}\draw(c)--(x)--(y)--(z);\draw(c)--(w)--(a);\eopfan{c}{180}{s-2}}
		\hfill\eoppanel{.31\linewidth}{$\mathcal R_1$(b): two leaves at $x$\\$s\ge3$}{\eopv{c}{0}{0}\eopv{x}{.6}{.35}\eopv{y}{1.1}{.6}\eopv{z}{1.15}{.15}\eopv{w}{.6}{-.4}\eopv{a}{1.15}{-.6}\draw(c)--(x)--(y);\draw(x)--(z);\draw(c)--(w)--(a);\eopfan{c}{180}{s-2}}
		\hfill\eoppanel{.31\linewidth}{$\mathcal R_1$(c): subdivisions\\$s\ge3$}{\eopv{c}{0}{0}\foreach \ang/\id in {60/a,0/b,-60/d}{\node[eop vertex](\id)at($(c)+(\ang:.55)$){};\node[eop vertex](\id e)at($(c)+(\ang:1.1)$){};\draw(c)--(\id)--(\id e);}\eopfan{c}{180}{s-3}}
		\caption{The six templates defining $\mathcal R_1$. In (a), a leaf of the initial star is identified with the specified vertex of $K_3$, $P_4$, or $K_{1,3}$. In (b), the two displayed branches are modified; in (c), three distinct star edges are subdivided once.}
		\label{fig:Families_graphs_R_class_part_1}
\end{figure}}

\newcommand{\EOPFigureRrest}{%
	\begin{figure}[tbp]\centering
		\eoppanel{.23\linewidth}{$\mathcal R_2$(a)}{\eoptriangle\eopfan{u}{90}{a\ge1}\eopv{a}{-.52}{-.95}\eope{v}{a}}
		\hfill\eoppanel{.23\linewidth}{$\mathcal R_3$(a), opposite}{\eopsquare\eopfan{u}{90}{a\ge1}\eopv{a}{0}{-1.1}\eope{w}{a}}
		\hfill\eoppanel{.23\linewidth}{$\mathcal R_3$(a), adjacent}{\eopsquare\eopfan{u}{90}{a\ge1}\eopv{a}{1.1}{0}\eope{z}{a}}
		\hfill\eoppanel{.23\linewidth}{$\mathcal R_2$(b)}{\eoptriangle\eopfan{u}{90}{a}\eopv{a}{.6}{.8}\eopv{b}{1.1}{1.1}\draw(u)--(a)--(b);}
		\par\medskip
		\eoppanel{.23\linewidth}{$\mathcal R_3$(b)}{\eopsquare\eopfan{u}{90}{a}\eopv{a}{.6}{.9}\eopv{b}{1.1}{1.2}\draw(u)--(a)--(b);}
		\hfill\eoppanel{.23\linewidth}{$\mathcal R_4$: $\ell=1$}{\eoptheta{1}\eopfan{u}{90}{a}}
		\hfill\eoppanel{.23\linewidth}{$\mathcal R_4$: $\ell=2$}{\eoptheta{2}\eopfan{u}{90}{a}}
		\hfill\eoppanel{.23\linewidth}{$\mathcal R_5$}{\foreach \ang/\id in {90/u,162/v,234/w,306/x,378/y}{\node[eop vertex](\id)at(\ang:.6){};}\draw(u)--(v)--(w)--(x)--(y)--(u);\eopfan{u}{90}{a}}
		\caption{Cycle-based templates. The first three panels add leaves at two cycle vertices; the next two add a path of length two at one cycle vertex. The two $\mathcal R_4$ templates are theta graphs with branch lengths $(\ell,2,2)$, where $\ell\in\{1,2\}$. Unspecified leaf counts may be zero.}
		\label{fig:Families_graphs_R_class_part_2}
\end{figure}}

\newcommand{\EOPFigureSone}{%
	\begin{figure}[tbp]\centering
		\eoppanel{.47\linewidth}{$\mathcal S_1$(a): $q=5$, $a,b\ge1$}{\eoppath{5}{4}\eopfan{v2}{90}{a}\eopfan{v4}{90}{b}\eopv{x}{1.44}{.55}\eope{v3}{x}}
		\hfill\eoppanel{.47\linewidth}{$\mathcal S_1$(b): $q\in\{6,7\}$}{\eoppath{6}{5}\eopfan{v2}{90}{a}\eopfan{v5}{90}{b}\eopv{x}{1.44}{.55}\eope{v3}{x}}
		\par\medskip
		\eoppanel{.47\linewidth}{$\mathcal S_1$(c): $q=7$}{\eoppath{7}{6}\eopfan{v2}{90}{a}\eopfan{v6}{90}{b}\eopv{x}{1.92}{.55}\eope{v4}{x}}
		\hfill\eoppanel{.47\linewidth}{$\mathcal S_1$(d): $q=8$}{\eoppath{8}{7}\eopfan{v2}{90}{a}\eopfan{v7}{90}{b}}
		\par\medskip
		\eoppanel{.65\linewidth}{$\mathcal S_1$(e): $q\in\{6,7,8\}$, $b\ge1$}{\eoppath{8}{7}\eopfan{v2}{90}{a}\eopfan{v6}{90}{b}\node[eop label,below=3pt]at(v2){$v_2$};\node[eop label,below=3pt]at(v6){$v_{q-2}$};}
		\caption{The path templates defining $\mathcal S_1$. The two support vertices have independent added-leaf counts. The drawings in (b) and (e) show one permitted path order; the ranges of $q$ are specified below the panels. The extra pendant edge in (b) may be placed symmetrically at either end.}
		\label{fig:Families_graphs_S_class_part_1}
\end{figure}}

\newcommand{\EOPFigureSrest}{%
	\begin{figure}[tbp]\centering
		\eoppanel{.30\linewidth}{$\mathcal S_2$}{\eoptriangle\eopv{x}{0}{1.0}\eopv{y}{0}{1.55}\draw(u)--(x)--(y);\eopfan{y}{90}{a}}
		\hfill\eoppanel{.30\linewidth}{$\mathcal S_3$(a)}{\eopsquare\eopv{x}{0}{1.1}\eopv{y}{0}{1.65}\draw(u)--(x)--(y);\eopfan{y}{90}{a}\eopfan{w}{270}{b}}
		\hfill\eoppanel{.30\linewidth}{$\mathcal S_3$(b)}{\eopsquare\eopv{x}{0}{1.1}\eope{u}{x}\eopfan{x}{90}{a}\eopfan{w}{270}{b}}
		\par\medskip
		\eoppanel{.46\linewidth}{$\mathcal S_4$: $(2,2,2)$}{\eoptriad{2}{2}{2}}
		\hfill\eoppanel{.46\linewidth}{$\mathcal S_4$: $(2,2,3)$}{\eoptriad{2}{2}{3}}
		\par\medskip
		\eoppanel{.46\linewidth}{$\mathcal S_4$: $(2,3,3)$}{\eoptriad{2}{3}{3}}
		\hfill\eoppanel{.46\linewidth}{$\mathcal S_4$: $(3,3,3)$}{\eoptriad{3}{3}{3}}
		\caption{The templates defining $\mathcal S_2$, $\mathcal S_3$, and $\mathcal S_4$. The last four panels list the three arm lengths of the underlying subdivided claw. The leaf counts $r_1,r_2,r_3$ are independent. In $\mathcal S_3$(b), $a\ge2$; other displayed leaf counts are nonnegative.}
		\label{fig:Families_graphs_S_class_part_2}
\end{figure}}

\newcommand{\EOPFigureThreeEdges}{%
	\begin{figure}[tbp]\centering
		\eoppanel{.30\linewidth}{$\mathcal G_1=K_{1,3}$}{\eopv{x}{0}{.8}\eopv{u}{-.9}{-.6}\eopv{y}{0}{-.6}\eopv{z}{.9}{-.6}\draw(x)--(u);\draw(x)--(y);\draw(x)--(z);\node[eop label]at(-.7,.15){$e_1$};\node[eop label]at(.17,-.2){$e_2$};\node[eop label]at(.65,.25){$e_3$};\node[eop label,above=3pt]at(x){$x$};\foreach\v in {u,y,z}{\node[eop label,below=3pt]at(\v){$\v$};}}
		\hfill\eoppanel{.30\linewidth}{$\mathcal G_2=P_3\cup K_2$}{\eopv{x}{0}{.8}\eopv{u}{-.65}{-.05}\eopv{y}{.65}{-.05}\eopv{z}{-.65}{-.9}\eopv{w}{.65}{-.9}\draw(x)--node[left=2pt,eop label]{$e_1$}(u);\draw(x)--node[right=2pt,eop label]{$e_2$}(y);\draw(z)--node[above=2pt,eop label]{$e_3$}(w);\foreach\v in {u,y,z,w}{\node[eop label,below=3pt]at(\v){$\v$};}\node[eop label,above=3pt]at(x){$x$};}
		\hfill\eoppanel{.30\linewidth}{$\mathcal G_3=P_4$}{\eopv{u}{.7}{.7}\eopv{x}{-.7}{.7}\eopv{y}{-.7}{-.7}\eopv{z}{.7}{-.7}\draw(u)--node[above=2pt,eop label]{$e_1$}(x);\draw(x)--node[left=2pt,eop label]{$e_2$}(y);\draw(y)--node[below=2pt,eop label]{$e_3$}(z);\node[eop label,above=3pt]at(u){$u$};\node[eop label,above=3pt]at(x){$x$};\node[eop label,below left=3pt]at(y){$y$};\node[eop label,below right=3pt]at(z){$z$};}
		\par\medskip
		\eoppanel{.40\linewidth}{$\mathcal G_4=C_3$}{\eopv{x}{0}{.8}\eopv{y}{-.75}{-.6}\eopv{z}{.75}{-.6}\draw(x)--node[left=2pt,eop label]{$e_1$}(y);\draw(y)--node[below=2pt,eop label]{$e_2$}(z);\draw(z)--node[right=2pt,eop label]{$e_3$}(x);\node[eop label,above=3pt]at(x){$x$};\node[eop label,below left=3pt]at(y){$y$};\node[eop label,below right=3pt]at(z){$z$};}
		\hfill\eoppanel{.40\linewidth}{$\mathcal G_5=3K_2$}{\foreach\y/\a/\b/\i in {.8/u/v/1,0/x/y/2,-.8/z/w/3}{\eopv{\a}{-.7}{\y}\eopv{\b}{.7}{\y}\draw(\a)--node[above=2pt,eop label]{$e_\i$}(\b);\node[eop label,below=3pt]at(\a){$\a$};\node[eop label,below=3pt]at(\b){$\b$};}}
		\caption{The five graphs formed by the three omitted edges. Only these edges and their endpoints are shown; other edges of $G$ are omitted.}
		\label{fig:Edge_induced_subgraph_by_the_three_edges}
\end{figure}}

\newcommand{\EOPFigureCertificates}{%
	\begin{figure}[tbp]\centering
		\eoppanel{.23\linewidth}{$3K_2$}{%
			\foreach\y/\a/\b in {.8/1/2,0/3/4,-.8/5/6}{\eopv{a\a}{-.65}{\y}\eopv{a\b}{.65}{\y}\draw(a\a)--(a\b);\node[eop label,left=3pt]at(a\a){$f_\a$};\node[eop label,right=3pt]at(a\b){$f_\b$};}}
		\hfill\eoppanel{.23\linewidth}{$K_3\cup K_2$}{%
			\eopv{a}{0}{1.05}\eopv{b}{-.6}{.2}\eopv{c}{.6}{.2}\eopv{d}{-.6}{-.8}\eopv{e}{.6}{-.8}\draw(a)--(b)--(c)--(a);\draw(d)--(e);\node[eop label,above=3pt]at(a){$f_1$};\node[eop label,left=3pt]at(b){$f_2$};\node[eop label,right=3pt]at(c){$f_3$};\node[eop label,left=3pt]at(d){$f_4$};\node[eop label,right=3pt]at(e){$f_5$};}
		\hfill\eoppanel{.23\linewidth}{$C_5$}{%
			\foreach\ang/\id in {90/1,162/2,234/3,306/4,378/5}{\node[eop vertex](a\id)at(\ang:.8){};\node[eop label]at(\ang:1.1){$f_\id$};}\draw(a1)--(a2)--(a3)--(a4)--(a5)--(a1);}
		\hfill\eoppanel{.23\linewidth}{$K_4$}{%
			\eopv{a}{-.65}{.65}\eopv{b}{.65}{.65}\eopv{c}{-.65}{-.65}\eopv{d}{.65}{-.65}\draw(a)--(b)--(d)--(c)--(a);\draw(a)--(d);\draw(b)--(c);\node[eop label,above=3pt]at(a){$f_1$};\node[eop label,above=3pt]at(b){$f_2$};\node[eop label,below=3pt]at(c){$f_3$};\node[eop label,below=3pt]at(d){$f_4$};}
		\caption{The four local certificates. A vertex labelled $f_i$ represents an edge of $G$; a drawn edge requires the corresponding pair to conflict. Additional conflicts are allowed, so these are subgraph patterns rather than induced-subgraph patterns. Each pattern requires at least three vertex deletions to destroy every conflict.}
		\label{fig:conflict-certificates}
\end{figure}}

\tikzset{eop edge/.style={line width=.45pt}}

\title{Edge open packing: structural characterizations, parameterized algorithms, and equality with induced matching}

\author[1]{arti@iitrpr.ac.in\footnote{arti@iitrpr.ac.in}}
	
	\affil[1]{Department of Mathematics\\
		
			Indian Institute of Technology Ropar\\
			
			 Rupnagar, 140001, Punjab} 
	
\author[2]{Kamal Santra \orcidlink{0009-0006-5997-1452} \footnote{kamal.7.2013@gmail.com, kamal.santra@iitg.ac.in}}

\affil[2]{Department of Mathematics\\
	
	Indian Institute of Technology Guwahati\\
	
	Guwahati, 781039, Assam, India}

\date{}

\begin{document}

	\maketitle
\begin{abstract}
We study graph structure through the edge open packing number $\rhooe(G)$. A packing selects edges whose endpoint-induced subgraph consists precisely of those edges arranged in star components. For connected graphs, we give a criterion for every positive value of $\rhooe(G)$. We also characterize $\rhooe(G)=m-3$, where $m=|E(G)|$, by a common construction from disjoint stars and four conflict patterns. The target-size decision problem is $\mathrm{W}[1]$-complete, while the deficit parameter admits a quadratic kernel and an $O(n+m+2^k(k+1)^4)$-time algorithm. Thus, every fixed deficit can be recognized in linear time. Finally, we extend a known equality construction for trees. Starting with disjoint subdivided stars, we add edges incident with their centres while keeping two original leaves per star pendant. The resulting graph has equal induced matching and edge open packing numbers, and the original pendant edges form its unique maximum induced matching. This class contains every finite graph as an induced subgraph.
\end{abstract}


\section{Introduction}\label{sec:introduction}
Choosing mutually separated edges gives an induced matching. Edge open packing allows several selected edges to meet at a centre: the leaves of each selected star must be pairwise nonadjacent, and different selected stars must have no edges between their endpoints. Chelladurai et al.~\cite{chelladurai2022edge} introduced the parameter $\rhooe(G)$ for the maximum number of edges selected under these rules. If $\nuI(G)$ is the induced matching number, then $\nuI(G)\leq\rhooe(G)$. For the star $K_{1,r}$, the two values are one and $r$, respectively, so their difference is unbounded.

We approach the parameter through three structural questions. Which graphs have a prescribed optimum $t$? What can be said when only $k$ of the $m$ edges must be excluded? When does permitting a shared centre give no improvement over an induced matching? The first question constrains every choice of $t+1$ edges, while the second asks how far the graph is from an induced star forest on its selected endpoints. We develop characterizations for these settings and examine the algorithms that follow from them.

\subsection{Related work}

Induced matching provides the natural starting point for this study. Its complexity was investigated by Stockmeyer and Vazirani~\cite{stockmeyer1982np}, and algorithms and hardness results for several graph classes were subsequently obtained in~\cite{cameron1989induced,lozin2002maximum,dabrowski2013new}. Replacing the selected matching by selected star components leads to the edge open packing problem.

The foundational work of Chelladurai et al.~\cite{chelladurai2022edge} gives bounds and characterizations when the optimum is one or two, or when its distance from $m$ is at most two. It also asks which graphs satisfy $\nuI(G)=\rhooe(G)$. Bre\v{s}ar and Samadi~\cite{brevsar2024edge} proved NP-completeness under several restrictions: a universal vertex, Eulerian bipartite structure, and planarity with maximum degree four. Their positive results include a linear-time tree algorithm and a characterization of equality in a bound involving minimum degree. Other structural results concern unicyclic graphs with optimum three~\cite{GayathriVenkataraman2025} and bounds for trees, split graphs, and unicyclic graphs~\cite{GayathriVenkataraman2026}.

For the equality question, Bre\v{s}ar et al.~\cite{BRESAR2025} completely described the trees. Their nontrivial construction starts with disjoint subdivided stars and connects them through edges meeting original centres; two original leaves per star must remain pendant. They also investigate five graph products: lexicographic, direct, Cartesian, strong, and rooted. Our equality result retains the leaf condition while allowing cycles and arbitrary induced cores. It yields an exact rooted-product formula when the attached graph is a subdivided star rooted at its centre, a special case of the rooted-product question in~\cite{BRESAR2025}.

Edge coloring gives a second connection. A color class in an injective edge coloring is an EOP, so minimizing the number of colors partitions the edge set into as few EOPs as possible~\cite{cardoso2019injective}. This covering parameter also occurs as induced star arboricity~\cite{AxenovichDorrRollinUeckerdt2019}. Algorithms and complexity for the coloring problem are studied in~\cite{foucaud2021complexity,PriyamvadaPanda2023}, and further bounds appear in~\cite{miao2022note}. Here the objective is the size of a single EOP.

Several algorithms exploit restrictions on the input graph. Santra's bounds for proper interval, block, and split graphs are, respectively, $O(n^2)$, $O(n+m)$, and $O(n^3)$~\cite{Santra2025Chordal}. Distance-hereditary and biconvex bipartite inputs are handled in polynomial time by Bhyravarapu et al.~\cite{BhyravarapuDasSantra2026}; their chordal-graph algorithm takes $2^{O(\omega)}\operatorname{poly}(n)$ time, where $\omega$ denotes clique number. For permutation and interval inputs with representations, the oriented star-conflict construction of Das and Santra gives a bound of $O(n^2+m^2)$~\cite{DasSantra2026Permutation}. Their separate algorithm for well-partitioned chordal graphs uses a partition-tree representation and takes $O(n^4)$ time.

Subsequent work of Das and Santra treats cactus and $P_4$-tidy graphs in $O(n+m)$ time, and cocomparability graphs in $O(n^2+m^2)$ time~\cite{DasSantra2026Cactus}. Their AT-free result relies on preservation of that class by the oriented construction and runs in $O(n^2+m^4)$ time~\cite{DasSantra2026ATFree}. The latter paper also gives separate $O(n^4)$ algorithms for chordal and convex bipartite graphs. These results concern unrestricted optimum sizes within specified graph classes. We study restrictions on the optimum itself and on the number of omitted edges.

For the parameterized analysis, we use the standard framework of~\cite{CyganEtAl2015}. The deficit is a Vertex Cover objective in an edge-conflict graph. The treewidth consequence of monadic second-order optimization~\cite{CourcelleMosbah1993}, already noted for EOP in~\cite{DasSantra2026ATFree}, is recalled with an explicit decomposition requirement.

\subsection{Our contributions}
We first give necessary and sufficient conditions for $2\leq\rhooe(G)\leq t$ for every integer $t\geq2$, yielding an exact-value characterization for each $t\geq3$. The proof treats induced matchings and larger induced stars in a common argument.

For the near-extremal value $m-3$, we give a common construction from at most three disjoint stars and three additional edges, together with four finite conflict certificates. The certificates are subgraph patterns in the edge-conflict graph. We also describe and verify the concrete star, cycle, path, and three-branch families. Every connected graph with a deficit of three has a vertex cover of at most four and at most two independent cycles; we describe the bicyclic case explicitly.

The target-size decision problem is $\mathrm{W}[1]$-complete, even on triangular graphs of diameter two with two universal vertices and with $\nuI(G)=\rhooe(G)$. This follows from a known clique-substitution identity. The deficit parameter $k$ admits a kernel with at most $3k(k+2)$ vertices and $O(k^2)$ edges, computable in $O(n+m)$ time; the vertex bound improves to $(2k+1)(k+2)$ for connected inputs. Vertex Cover branching on this kernel gives a deterministic $O(n+m+2^k(k+1)^4)$-time algorithm and reconstructs an optimum when the deficit is at most $k$. Thus, every fixed deficit can be recognized in linear time. The kernel also gives an effectively computable finite-core construction for every fixed deficit. We give a further kernel for the combined parameter $t+\Delta(G)$ and record the treewidth consequence of monadic second-order optimization.

Finally, we extend the spider construction of Bre\v{s}ar et al.~\cite{BRESAR2025}. Start with disjoint subdivided stars and add arbitrary edges incident with their original centres, keeping at least two original leaves of each star pendant. Both parameters equal the total number of arms, and the original pendant edges form the unique maximum induced matching. The construction permits cycles and arbitrary induced graph cores.

Section~\ref{sec:preliminaries} fixes notation. Section~\ref{EOP_sec_2} treats prescribed EOP numbers. Section~\ref{EOP_sec_3} gives the common characterization for deficit three and all concrete constructions with their proofs. Section~\ref{sec:parameterized} presents the parameterized results. Section~\ref{sec:equal-parameters} proves the equality extension, and Section~\ref{EOP_sec_4} concludes the paper.
\section{Preliminaries}\label{sec:preliminaries}
We consider finite simple graphs and put $n=|V(G)|$, $m=|E(G)|$. The structural characterizations assume connectedness and at least one edge unless stated otherwise. Our notation is $P_q$ for a path of order $q$ and $C_q$ for a cycle of order $q$. Given $X\subseteq V(G)$, the notation $G[X]$ specifies the induced subgraph. For $D\subseteq E(G)$, write $V_D$ for all endpoints of edges in $D$. We abbreviate $G[V_D]$ to $G[D]$, whereas $G\langle D\rangle$ has precisely the edges in $D$. The symbols $N_G(v)$, $\operatorname{diam}(G)$, and $\tau(G)$ denote open neighbourhood, diameter, and vertex-cover number. A vertex adjacent to a leaf is a support vertex. Further terminology follows West~\cite{west2001introduction}.

The EOP condition can be stated as follows: $E(G[V_D])=D$, and every component of this induced subgraph is a nontrivial star. We may choose either endpoint as centre when a component is $K_2$. Values of $\nuI$ and $\rhooe$ sum over components of $G$, with zero contribution from an isolated vertex.

To describe incompatible choices, call two distinct edges $e,f$ \emph{conflicting} if some $g\notin\{e,f\}$ has one endpoint incident with $e$ and its other endpoint incident with $f$. Taken together, $e,f,g$ form either a three-edge path or a triangle. Thus, selecting an EOP is equivalent to selecting edges without a conflicting pair. We use $\alpha(G)$ and $\Delta(G)$ for independence number and maximum degree, and write $\eta(G)=m-\rhooe(G)$ for the deficit.

Parameterized terminology follows~\cite{CyganEtAl2015}. For a parameter $p$, an FPT running time has the form $f(p)(n+m)^c$ for a constant $c$ independent of $p$; an XP running time may have $p$ in the exponent. Kernelization replaces an input by an equivalent instance in polynomial time, with the output size controlled solely by the parameter. Our kernels retain the EOP formulation: the reduced instance consists of a graph and a parameter.

\section{Prescribed edge open packing numbers}\label{EOP_sec_2}
For a connected graph with at least one edge, $\rhooe(G)=1$ if and only if $G$ is complete, and $\rhooe(G)=m$ if and only if $G$ is a star~\cite{chelladurai2022edge}. We also use the bound $\rhooe(G)\geq\ceil{\operatorname{diam}(G)/2}$ from~\cite{chelladurai2022edge}. The next theorem includes the known characterization for $\rhooe(G)=2$ when $t=2$.

\begin{theorem}\label{EOP_any_char_theorem}
	Let $G$ be connected and let $t\geq2$. Then $2\leq\rhooe(G)\leq t$ if and only if the following conditions hold.
	\begin{enumerate}[label=\textup{(C\arabic*)},leftmargin=*]
		\item $2\leq\operatorname{diam}(G)\leq2t$.
		\item $G$ has no induced $K_{1,t+1}$.
		\item For every induced matching $M$ of size $t$, each vertex outside $V_M$ has at least two neighbours in $V_M$.
		\item For every EOP set $F$ of size $t$ whose induced star forest has between two and $t-1$ components, each vertex $z\notin V_F$ adjacent to a centre $u$ of a component has another neighbour in $V_F\setminus\{u\}$. This is required for either choice of centre in a $K_2$ component. For $t=2$, this condition is empty.
	\end{enumerate}
\end{theorem}

\begin{proof}
	Assume that $2\leq\rho_e^o(G)\leq t$. The diameter bound and the characterization of graphs with EOP number one give $(C_1)$. An induced $K_{1,t+1}$ supplies an EOP set of size $t+1$, so $(C_2)$ follows.
	
	To prove $(C_3)$, let $M=\{x_iy_i:1\leq i\leq t\}$ be an induced matching, put $S=V(G[M])$, and assume that $z\notin S$ has at most one neighbour in $S$. If its unique neighbour is $x$, then $M\cup\{zx\}$ is an EOP set of size $t+1$. Otherwise, $z$ has no neighbour in $S$. Since $G$ is connected, choose $w\in N_G(z)$. Put $I=\{i:N_G(w)\cap\{x_i,y_i\}\neq\emptyset\}$ and, for each $i\in I$, choose $a_i\in N_G(w)\cap\{x_i,y_i\}$. The set $(M\setminus\{x_iy_i:i\in I\})\cup\{wz\}\cup\{wa_i:i\in I\}$ is an EOP set of size $t+1$: the new edges form an induced star, and its vertices have no neighbours among the endpoints of the retained matching edges. Both cases contradict $\rho_e^o(G)\leq t$.
	
	For $(C_4)$, put $S=V_F$. If $N_G(z)\cap S=\{u\}$, then $F\cup\{zu\}$ is an EOP set of size $t+1$, again a contradiction.
	
	Conversely, assume $(C_1)$--$(C_4)$ and suppose that an EOP set $D$ of size $t+1$ exists. The graph $G[D]$ is an induced star forest; let $k$ be its number of components. If $k=1$, then $G[D]\cong K_{1,t+1}$, contrary to $(C_2)$. If $k=t+1$, choosing one edge from each of $t$ components gives an induced matching $M$ whose endpoints have no neighbour in the remaining component, contrary to $(C_3)$.
	
	It remains that $2\leq k\leq t$. Some component of $G[D]$ has at least two edges. Delete one of its leaf edges $zu$ and put $F=D\setminus\{zu\}$. Then $|F|=t$, $G[F]$ still has $k$ components, $z\notin V(G[F])$, and $N_G(z)\cap V(G[F])=\{u\}$. If $k=t$, then $F$ is an induced matching and contradicts $(C_3)$. If $2\leq k\leq t-1$, choose $u$ as the centre of its component of $G[F]$; this contradicts $(C_4)$. Thus $\rho_e^o(G)\leq t$. Finally, $(C_1)$ implies that $G$ is not complete, and hence $\rho_e^o(G)\geq2$.
\end{proof}

For $q\geq2$, let $\mathcal F_q$ denote the conjunction of the conditions of Theorem~\ref{EOP_any_char_theorem} with $t=q$.
\begin{corollary}\label{cor:exact-t}
	Let $G$ be connected and let $t\geq3$. Then $\rhooe(G)=t$ if and only if $G$ satisfies $\mathcal F_t$ and fails at least one condition of $\mathcal F_{t-1}$.
\end{corollary}
\begin{proof}
	The two sets of conditions mean $2\leq\rhooe(G)\leq t$ and $\rhooe(G)\notin[2,t-1]$, respectively.
\end{proof}
The corollary requires no separate exclusion of stars. Section~\ref{sec:parameterized} addresses its algorithmic consequences when $t$ is part of the input.

\section{A unified characterization of graphs with \texorpdfstring{$\rhooe(G)=m-3$}{EOP number m-3}}\label{EOP_sec_3}
We first recall the near-extremal context. The connected graphs with deficit $m-\rhooe(G)$ equal to zero are precisely the stars. For deficit one, they are obtained by subdividing one edge of $K_{1,m-1}$ once, where $m\geq3$. Chelladurai et al.~\cite{chelladurai2022edge} also characterized deficit two. In place of a separate case list for deficit three, we use a common construction and four finite certificates. The construction describes the underlying graph; the certificate distinguishes the deficit three from the smaller deficits.

\subsection{The common construction}
\begin{definition}\label{def:extension}
	A \emph{$k$-edge star extension} is a connected graph obtained by taking a disjoint union of nontrivial stars, adding new isolated vertices, and then adding exactly $k$ edges, each of which has at least one new endpoint. No other edges or vertices are added.
\end{definition}
The new endpoints are new vertices, rather than newly designated vertices of a star. In particular, added edges cannot join two vertices of the original stars. A new vertex may be joined to several star vertices or to other new vertices. All graphs in the construction are simple.

\begin{lemma}\label{lem:star-completion}
	Let $G$ be connected with $m$ edges, and let $1\leq k<m$. Then $\rhooe(G)\geq m-k$ if and only if $G$ is a $k$-edge star extension. If there are $q$ original stars and $r$ new vertices, then $1\leq r\leq k$, $q+r\leq k+1$, and $q\leq k$.
\end{lemma}
\begin{proof}
	Choose an EOP set $D$ of size $m-k$. Its endpoints induce the original disjoint stars. Every edge outside $D$ has an endpoint outside $V_D$, because $G[D]$ has exactly the edges $D$. Thus $G$ has the required construction, with the vertices outside $V_D$ as the new vertices. Since $k\geq1$, at least one such vertex exists. Contracting each original star gives a connected multigraph on $q+r$ vertices with $k$ edges. Hence $q+r\leq k+1$, which implies the stated bounds. Conversely, the original star edges form an EOP of size $m-k$, since no added edge has both endpoints in the original vertex set.
\end{proof}
For $k=3$, there are at most three stars and at most $4-q$ new vertices. Thus, the entire construction has a skeleton on at most four vertices after the stars are contracted. The three added edges, before contraction, form one of the five graphs in Figure~\ref{fig:Edge_induced_subgraph_by_the_three_edges}; no separate case analysis of those shapes is needed.
\EOPFigureThreeEdges

\begin{corollary}\label{cor:bounded-cover}
	If $G$ is connected and $\rhooe(G)=m-k$ with $1\leq k<m$, then $\tau(G)\leq k+1$. In particular, $\rhooe(G)=m-3$ implies $\tau(G)\leq4$.
\end{corollary}
\begin{proof}
	The original star centres and the new vertices in Lemma~\ref{lem:star-completion} cover every edge. Their total number is at most $k+1$.
\end{proof}

\subsection{Four finite conflict certificates}
Define the \emph{edge-conflict graph} $J(G)$ to have vertex set $E(G)$, with two vertices adjacent precisely when their corresponding edges conflict in $G$. An EOP in $G$ is an independent set in $J(G)$. Consequently, $m-\rhooe(G)=\tau(J(G))$.

\begin{lemma}\label{lem:cover-two}
	A graph $H$ has $\tau(H)\geq3$ if and only if it contains, as a subgraph, at least one of $3K_2$, $K_3\cup K_2$, $C_5$, and $K_4$. The subgraph is not required to be induced.
\end{lemma}
\begin{proof}
	Each of the four graphs has a vertex-cover number of three, so its presence implies $\tau(H)\geq3$. For the converse, assume that $H$ contains none of them. In particular, its matching number is at most two.
	
	If $H$ is triangle-free, it has no odd cycle: a shortest odd cycle would have length five, or length at least seven and hence contain three disjoint edges. Thus $H$ is bipartite. By K\H{o}nig's theorem, its vertex-cover number equals its matching number and is at most two.
	
	Otherwise, fix a triangle $T=abc$. The vertices outside $T$ are independent, since an edge between two of them would give $K_3\cup K_2$. If one triangle vertex has no neighbour outside $T$, the other two vertices cover every edge. Assume, therefore, that each triangle vertex has an external neighbour. No external vertex can see all three triangle vertices, since that would give a $K_4$. If an external vertex $x$ sees exactly two, say $a,b$, choose an external neighbour $y$ of $c$. Then $y\neq x$, and the triangle $abx$ and edge $cy$ give $K_3\cup K_2$. Hence, every external vertex has at most one neighbour in $T$. External neighbours chosen for $a,b,c$ are therefore distinct, giving three disjoint edges. This final contradiction proves $\tau(H)\leq2$.
\end{proof}

Figure~\ref{fig:conflict-certificates} shows these patterns. Table~\ref{tab:certificates} expresses these four possibilities directly in terms of edges of $G$. The edges listed in a row must be distinct. Other conflicts among them are permitted.
\begin{table}[htbp]
	\centering
	\caption{Four certificates that the EOP deficit is at least three.}
	\label{tab:certificates}
	\begin{tabular}{@{}l>{\raggedright\arraybackslash}p{.73\linewidth}@{}}
		\toprule
		Pattern in $J(G)$ & Required conflicts between edges of $G$\\
		\midrule
		$3K_2$ & Six edges $e_1,\ldots,e_6$ with conflicting pairs $(e_1,e_2)$, $(e_3,e_4)$, and $(e_5,e_6)$.\\[3pt]
		$K_3\cup K_2$ & Five edges $e_1,\ldots,e_5$; the first three are pairwise conflicting, and $e_4,e_5$ conflict.\\[3pt]
		$C_5$ & Five edges $e_1,\ldots,e_5$, with $e_i$ conflicting with $e_{i+1}$, where indices are read modulo five.\\[3pt]
		$K_4$ & Four pairwise conflicting edges.\\
		\bottomrule
	\end{tabular}
\end{table}

\EOPFigureCertificates

\begin{definition}\label{def:four-families}
	For $F\in\{3K_2,K_3\cup K_2,C_5,K_4\}$, let $\mathcal B_F$ consist of the graphs constructed as follows:
	\begin{enumerate}[label=\textup{(\roman*)},leftmargin=*]
		\item take $q$ disjoint stars $K_{1,s_1},\ldots,K_{1,s_q}$, where $1\leq q\leq3$ and each $s_i\geq1$;
		\item add $r$ new vertices, where $1\leq r\leq4-q$;
		\item add exactly three edges, each with a new endpoint, so that the graph becomes connected and contains the conflict pattern $F$ in Table~\ref{tab:certificates}.
	\end{enumerate}
\end{definition}
These are four families with one common construction rule. The parameters $s_i$ are independent. The final condition involves at most six specified edges and their conflicts; it contains no reference to the unknown optimum $\rhooe(G)$. The families can overlap, because several certificates may occur in the same graph.

\begin{theorem}\label{EOP_char_theorem_3}
	Let $G$ be connected with $m\geq4$ edges. Then $\rhooe(G)=m-3$ if and only if $G\in\bigcup_{F\in\{3K_2,K_3\cup K_2,C_5,K_4\}}\mathcal B_F$.
\end{theorem}
\begin{proof}
	If $\rhooe(G)=m-3$, Lemma~\ref{lem:star-completion} gives the common construction with three added edges. Also, $\tau(J(G))=3$, so Lemma~\ref{lem:cover-two} supplies one of the four certificates. Thus $G$ belongs to the stated union.
	
	Conversely, the original star edges in a construction from Definition~\ref{def:four-families} form an EOP of size $m-3$. The certificate and Lemma~\ref{lem:cover-two} imply $\tau(J(G))\geq3$, so $\rhooe(G)\leq m-3$. Equality follows.
\end{proof}
The star construction bounds the deficit from above by three; a conflict certificate bounds it from below. The characterization therefore covers all placements of the three added edges, and the four families may overlap.

Let $\Theta_{\ell,2,2}$ denote three internally vertex-disjoint paths with common endpoints and lengths $\ell,2,2$.
\begin{corollary}\label{cor:m3-cyclomatic}
	A connected graph $G$ satisfying $\rho_e^o(G)=m-3$ is a tree, a
	unicyclic graph, or a bicyclic graph. In the bicyclic case, the graphs
	are exactly those obtained from $\Theta_{\ell,2,2}$, with
	$\ell\in\{1,2\}$, by adding arbitrarily many leaves at one of the
	two common endpoints of its three paths.
\end{corollary}
\begin{proof}
	In the common construction, let $q$ be the number of original stars
	and let $w$ be the number of new vertices. Since the stars contain
	$m-3$ edges, they have $m-3+q$ vertices. Thus
	$n=m-3+q+w$, so the cyclomatic number of $G$ is
	$m-n+1=4-q-w\leq2$.
	
	If this number is two, then $q=w=1$. Let $c$ be the centre of the
	original star and let $z$ be the unique new vertex. Its three added
	edges join $z$ to three distinct star vertices. If one is $c$, the
	other two are leaves and the resulting core is $K_4-e$, equivalently
	$\Theta_{1,2,2}$. Otherwise all three are leaves and the core is
	$K_{2,3}$, equivalently $\Theta_{2,2,2}$. Every remaining original
	star vertex is a leaf at $c$.
	
	Conversely, both types have the common star construction. In the
	first type, write the degree-three core vertices as $c,z$ and the
	other vertices as $x,y$. The edges $cz,cx,zx$ are pairwise
	conflicting, while $cy,zy$ form another conflicting pair, giving
	the $K_3\cup K_2$ certificate. In the second type, write the common
	neighbours of $c,z$ as $x,y,t$. The three pairs
	$\{cx,zy\}$, $\{cy,zt\}$, and $\{ct,zx\}$ are conflicting pairs,
	giving the $3K_2$ certificate. Additional leaves preserve these
	certificates. The characterization theorem gives the result.
\end{proof}

\begin{corollary}\label{cor:recognition}
	The equality $\rhooe(G)=m-3$ can be recognized in $O(n+m)$ time. When it holds, a maximum EOP can be constructed within the same bound.
\end{corollary}
\begin{proof}
	Apply Corollary~\ref{cor:deficit-algorithm} with the fixed deficit bound $k=3$. The kernel and reconstruction are given in Section~\ref{sec:parameterized}.
\end{proof}

\subsection{A triangle with distinct pendant structures}\label{subsec:omission}
Let $Q$ be the triangle $012$ with an added leaf $4$ at vertex $1$ and a path $035$ of length two at vertex $0$. It has six edges, and $\{12,14,35\}$ is an EOP. Moreover, $(01,35)$, $(02,12)$, and $(03,14)$ are three disjoint conflicting pairs, witnessed by the edges $03$, $01$, and $01$, respectively. Hence $J(Q)$ contains $3K_2$ and $\rhooe(Q)=3=|E(Q)|-3$.

The two pendant structures occur at distinct triangle vertices. This example illustrates why the common construction permits independent attachment choices. The named constructions in Section~\ref{sec:families} are illustrative subclasses; the full characterization is Theorem~\ref{EOP_char_theorem_3}.

\begin{figure}[htbp]
	\centering
	\begin{tikzpicture}[x=1cm,y=1cm,every node/.style={font=\small}]
		\node[eop vertex,label=below:$0$] (a) at (0,0) {};
		\node[eop vertex,label=below:$1$] (b) at (1.4,0) {};
		\node[eop vertex,label=above:$2$] (c) at (.7,1.1) {};
		\node[eop vertex,label=below:$3$] (d) at (-1.2,0) {};
		\node[eop vertex,label=below:$5$] (f) at (-2.4,0) {};
		\node[eop vertex,label=below:$4$] (e) at (2.6,0) {};
		\draw[eop edge] (a)--(b) (a)--(c) (a)--(d);
		\draw[eop edge,line width=1.1pt] (b)--(c) (b)--(e) (d)--(f);
	\end{tikzpicture}
	\caption{The graph $Q$. The three thick edges form a maximum EOP.}
	\label{fig:omitted-graph}
\end{figure}

\subsection{Concrete graph families and their verification}\label{sec:families}
We now describe concrete subclasses of Theorem~\ref{EOP_char_theorem_3} and verify them explicitly. They need not be disjoint or exhaustive. Every numerical parameter is independent unless a relation is stated.

\subsubsection{The recalled deficit-two families}
Write $P_q=v_1\cdots v_q$. The families $\mathcal A_1,\mathcal A_2,\mathcal A_3$ from~\cite{chelladurai2022edge} use $q=5,6,7$, respectively: add independently $a,b\geq0$ leaves at $v_2,v_{q-1}$. The families $\mathcal A_4,\mathcal A_5$ are obtained from $C_3,C_4$ by adding $a\geq0$ leaves at one vertex. For $s\geq3$, $\mathcal A_6$ is obtained from $K_{1,s}$ by subdividing two distinct edges once each, and $\mathcal A_7$ by adding two leaves at one original leaf. Figure~\ref{fig:Families_graphs_A_class} gives a common graphical convention.
\EOPFigureA

\subsubsection{Star, cycle, path, and three-branch constructions}
We use the following conventions throughout the constructions. To
\emph{add $a$ leaves at $x$} means to add $a$ new vertices, each adjacent
only to $x$. To \emph{add a path of length $\ell$ at $x$} means to take a
disjoint copy of $P_{\ell+1}$ and identify one endpoint with $x$. Thus, a
path of length two adds two vertices and two edges. All vertices added
in different operations are distinct, and leaf multiplicities at
different vertices may be chosen independently. Zero leaves means that
no vertex or edge is added.

\paragraph{Star constructions}
The family $\mathcal R_1$ consists of the following graphs obtained from
a star $K_{1,s}$.
\begin{enumerate}[label=(\alph*),leftmargin=*]
	\item Identify a leaf of $K_{1,s}$ with the specified vertex of one of
	the three graphs in Table~\ref{tab:R1-roots}.
	\item For $s\geq3$, choose two distinct leaves $x,y$ of $K_{1,s}$. Add a
	leaf at $y$ and, at $x$, either add two leaves or add a path of length
	two.
	\item For $s\geq3$, subdivide three distinct edges of $K_{1,s}$, each
	exactly once.
\end{enumerate}
\begin{table}[htbp]
	\centering
	\caption{The rooted identifications defining $\mathcal R_1$(a).
		The lower bound on $s$ is part of each construction.}
	\label{tab:R1-roots}
	\begin{tabular}{lll}
		\hline
		Graph & Vertex identified with the star leaf & Required $s$\\
		\hline
		$K_3$ & any vertex & $s\geq2$\\
		$P_4$ & a support vertex & $s\geq3$\\
		$K_{1,3}$ & its centre & $s\geq4$\\
		\hline
	\end{tabular}
\end{table}
Every construction of $\mathcal R_1$ increases the number of edges of the original star by exactly three. The bounds in Table~\ref{tab:R1-roots} exclude the
smaller stars whose resulting graphs have a deficit at most two; they do
not permit an unmodified star as an additional case.

\EOPFigureRone

\paragraph{Cycle constructions}
For $k\in\{3,4\}$, let $\mathcal C_k$ consist of the graphs obtained
from $C_k$ by one of the following two operations:
\begin{enumerate}[label=(\alph*),leftmargin=*]
	\item choose distinct cycle vertices $u,v$, add $a\geq1$ leaves at $u$,
	and add one leaf at $v$;
	\item choose a cycle vertex $u$, add $a\geq0$ leaves at $u$, and add a
	path of length two at $u$.
\end{enumerate}
We set $\mathcal R_2=\mathcal C_3$ and
$\mathcal R_3=\mathcal C_4$. For $k=4$, the vertices $u,v$ in (a) may
be adjacent or opposite; these are the two possible relative positions on the cycle.

Let $\Theta_{\ell,2,2}$ be the graph formed by three internally
vertex-disjoint paths with common endpoints $u,v$ and lengths
$\ell,2,2$, respectively. The family $\mathcal R_4$ consists of the
graphs obtained from $\Theta_{\ell,2,2}$, where $\ell\in\{1,2\}$, by
adding $a\geq0$ leaves at $u$. Equivalently, start with $K_4-e$,
subdivide its edge joining the two degree-three vertices zero or one
times, and add leaves at one of those vertices. 

The family $\mathcal R_5$ consists of the graphs obtained from $C_5$ by
adding $a\geq0$ leaves at one cycle vertex.

\EOPFigureRrest

\paragraph{Path constructions}
For a path $P_q=v_1\cdots v_q$, the family $\mathcal S_1$ is defined by
Table~\ref{tab:S1-leaves}. A row specifies all leaves added to the path;
no leaves are added at unlisted vertices. Path reversal gives the
equivalent choices at the other end. The parameters $a,b$ in each row
are independent.
\begin{table}[htbp]
	\centering
	\caption{A common description of the subfamilies of $\mathcal S_1$.}
	\label{tab:S1-leaves}
	\begin{tabular}{llll}
		\hline
		Subfamily & $q$ & Vertices and numbers of added leaves & Restrictions\\
		\hline
		(a) & $5$ & $v_2:a,\ v_4:b,\ v_3:1$ & $a,b\geq1$\\
		(b) & $6,7$ & $v_2:a,\ v_{q-1}:b,\ v_3:1$ & $a,b\geq0$\\
		(c) & $7$ & $v_2:a,\ v_6:b,\ v_4:1$ & $a,b\geq0$\\
		(d) & $8$ & $v_2:a,\ v_7:b$ & $a,b\geq0$\\
		(e) & $6,7,8$ & $v_2:a,\ v_{q-2}:b$ & $a\geq0,\ b\geq1$\\
		\hline
	\end{tabular}
\end{table}

\EOPFigureSone

\paragraph{A cycle with a rooted path or star}
The family $\mathcal S_2$ consists of the graphs obtained from $C_3$ by
choosing a vertex $u$, adding a path $uxy$ of length two at $u$, and
adding $a\geq0$ leaves at $y$.

For $\mathcal S_3$, start with $C_4$ and choose opposite vertices $u,v$.
Apply one of the following operations:
\begin{enumerate}[label=(\alph*),leftmargin=*]
	\item add a path $uxy$ of length two at $u$, add $a\geq0$ leaves at $y$,
	and add $b\geq0$ leaves at $v$;
	\item add a new vertex $w$ adjacent to $v$, add $a\geq2$ leaves at $w$,
	and add $b\geq0$ leaves at $u$.
\end{enumerate}
In (b), the new edge and leaves form a star with a cycle vertex as one
of its leaves; no identification or extra connecting edge is implicit.

\paragraph{Three branches}
The family $\mathcal S_4$ consists of the graphs obtained by identifying
one endpoint of three disjoint paths, whose lengths are independently
chosen from $\{2,3\}$, into a single vertex. At the neighbour of the
remaining endpoint of each path, add an arbitrary number of leaves.
The three leaf multiplicities are independent. Up to permuting the
branches, the four possibilities correspond to zero, one, two, or
three paths of length three.

The families need not be disjoint. In particular, a parameter value
zero may produce a graph also described by another family; this causes
no ambiguity because a graph may have more than one construction.

\EOPFigureSrest

\subsubsection{A combined cycle construction}

For $k\in\{3,4\}$ choose vertices $u,v$ of $C_k$, distinct when $k=3$
and opposite when $k=4$. Add any number $b\geq0$ of leaves at $v$.
Take a disjoint star $K_{1,r}$ with $r\geq1$ and join $u$ by a new
edge to a specified vertex of the star, which may be its centre or
a leaf. Every resulting graph has edge open packing number $m-3$.
The triangle version includes the graph in Section~\ref{subsec:omission}.
Indeed, the edges not incident with $u$ form the original star
together with a star centred at $v$, and hence give an edge open
packing of size $m-3$. To certify the reverse inequality, let $t$ be
the vertex of the disjoint star joined to $u$, and choose an original
star edge $ts$. For a triangle $uvw$, the edges $uv,uw,ts$ are
pairwise conflicting, and $vw,ut$ form a further conflicting pair.
For a four-cycle $upvqu$, the three pairs
$\{up,vq\}$, $\{pv,ut\}$, and $\{qu,ts\}$ are conflicting pairs.
Thus, one of the four certificates is present in either case.

\subsection{Verification of the named constructions}\label{subsec:named-certificates}
The following verification applies to all named constructions and explains the lower bounds on their parameters.

\begin{lemma}\label{lem:leaf-extension}
	Let $H$ have $h$ edges and an EOP set $D$ of size $h-k$. Choose one centre for each component of $H[D]$, and let $U$ be a subset of the chosen centres. If $J(H)$ has vertex-cover number at least $k$, then adding any number of new leaves independently at the vertices of $U$ produces a graph $G$ satisfying $\rhooe(G)=|E(G)|-k$.
\end{lemma}
\begin{proof}
	Extend $D$ by all new pendant edges. Its endpoints still induce a star forest, so the extended EOP has size $|E(G)|-k$. Every conflict between two old edges remains present, and hence $J(H)$ is a subgraph of $J(G)$. Thus $\tau(J(G))\geq\tau(J(H))\geq k$. The identity $\rhooe(G)=|E(G)|-\tau(J(G))$ gives the reverse inequality.
\end{proof}

\begin{proposition}\label{prop:named-families}
	Every graph in $\mathcal R_1,\ldots,\mathcal R_5$ or $\mathcal S_1,\ldots,\mathcal S_4$ has edge open packing number $m-3$.
\end{proposition}
\begin{proof}
	For each construction take its \emph{core} by setting every leaf multiplicity to its stated minimum. The full graph is obtained from its core by independently adding leaves at the vertices prescribed in its definition. Tables~\ref{tab:named-core-certificates} and~\ref{tab:path-core-certificates} provide a set $X$ of three core edges such that $D=E(H)\setminus X$ is an EOP, together with a conflict certificate from Lemma~\ref{lem:cover-two}. Each vertex at which further leaves may be added is a centre of a component of $H[D]$. Lemma~\ref{lem:leaf-extension}, with $k=3$, therefore applies.
	
	Here is the labelling used in Table~\ref{tab:named-core-certificates}. In a star construction, $c$ is the original centre, $x,y$ are the chosen original leaves, and $\ell_i$ are other original leaves. For $\mathcal R_1(a)$, the three identified graphs have edge sets $\{xa,ab,bx\}$, $\{xa,xb,bd\}$, and $\{xa,xb,xd\}$, respectively. In $\mathcal R_1(b)$ the extra edges are either $\{xa,xb,yd\}$ or $\{xa,ab,yd\}$, and $\ell$ is an unchosen original leaf. For $\mathcal R_1(c)$ the three subdivided branches are $cx_i y_i$, where $1\leq i\leq3$.
	
	The triangle and four-cycle are labelled $uvzu$ and $upvqu$, respectively. In $\mathcal R_2(a)$ the core pendant edges are $ua,vb$; in the adjacent and opposite versions of $\mathcal R_3(a)$ they are $ua,pb$ and $ua,vb$, respectively. In the path versions of these cycle constructions the added path is $uxy$. The labels for $\mathcal S_2$ and $\mathcal S_3$ agree with their definitions; in $\mathcal S_3(b)$, the two required new leaves at $w$ are $a,b$.
	
	In the certificate column, $[e,f]$ means that the two edges conflict. Three brackets give $3K_2$; $\langle e,f,g\rangle;[h,j]$ means that the first three edges are pairwise conflicting and the last pair conflicts, giving $K_3\cup K_2$. All edges appearing in one certificate are distinct. Thus, every row is directly verifiable from its core, without a comparison with the previously characterized extremal families.
	
	\begingroup
	\small
	\setlength{\tabcolsep}{4pt}
	\renewcommand{\arraystretch}{1.30}
	\begin{longtable}{@{}>{\raggedright\arraybackslash}p{.20\linewidth}>{\raggedright\arraybackslash}p{.25\linewidth}>{\raggedright\arraybackslash}p{.49\linewidth}@{}}
		\caption{Deletion triples and conflict certificates for the star and cycle constructions.}\label{tab:named-core-certificates}\\
		\toprule
		Construction & Deleted edges $X$ & Conflict certificate\\
		\midrule
		\endfirsthead
		\toprule
		Construction & Deleted edges $X$ & Conflict certificate\\
		\midrule
		\endhead
		$\mathcal R_1(a)$, $K_3$ & $xa,ab,bx$ & $\langle c\ell_1,xa,xb\rangle;[cx,ab]$\\
		$\mathcal R_1(a)$, $P_4$ & $xa,xb,bd$ & $[c\ell_1,xb];[c\ell_2,xa];[cx,bd]$\\
		$\mathcal R_1(a)$, $K_{1,3}$ & $xa,xb,xd$ & $[c\ell_1,xa];[c\ell_2,xb];[c\ell_3,xd]$\\
		$\mathcal R_1(b)$, two leaves at $x$ & $xa,xb,yd$ & $[cy,xb];[c\ell,xa];[cx,yd]$\\
		$\mathcal R_1(b)$,\newline path at $x$ & $xa,ab,yd$ & $[c\ell,yd];[cx,ab];[cy,xa]$\\
		$\mathcal R_1(c)$ & $x_1y_1,x_2y_2,x_3y_3$ & $[cx_1,x_2y_2];[cx_2,x_3y_3];[cx_3,x_1y_1]$\\
		$\mathcal R_2(a)$ & $uv,vz,vb$ & $\langle uv,vz,zu\rangle;[ua,vb]$\\
		$\mathcal R_2(b)$ & $uv,vz,xy$ & $\langle uv,uz,xy\rangle;[vz,ux]$\\
		$\mathcal R_3(a)$, adjacent roots & $pv,vq,pb$ & $[pv,ua];[up,vq];[uq,pb]$\\
		$\mathcal R_3(a)$, opposite roots & $pv,vq,vb$ & $[pv,ua];[up,vq];[uq,vb]$\\
		$\mathcal R_3(b)$ & $pv,vq,xy$ & $[pv,ux];[up,vq];[uq,xy]$\\
		$\mathcal S_2$ & $uv,uz,ux$ & $\langle uv,uz,xy\rangle;[vz,ux]$\\
		$\mathcal S_3(a)$ & $up,uq,ux$ & $[pv,ux];[up,vq];[uq,xy]$\\
		$\mathcal S_3(b)$ & $pv,vq,vw$ & $[uq,vw];[vq,wb];[pv,wa]$\\
		\bottomrule
	\end{longtable}
	\endgroup
	
	For $\mathcal S_1$, put $e_i=v_iv_{i+1}$. In row (a), let $a_0,b_0$ be the required pendant edges at $v_2,v_4$, and put $h=v_3z$. In rows (b) and (c), let $h$ be the one extra pendant edge specified in the definition. In row (e), let $h$ be one required pendant edge at $v_{q-2}$. Every certificate in Table~\ref{tab:path-core-certificates} consists of three distinct conflicting pairs.
	
	\begin{table}[htbp]
		\centering
		\small
		\setlength{\tabcolsep}{4pt}
		\renewcommand{\arraystretch}{1.30}
		\caption{Deletion triples and conflict certificates for $\mathcal S_1$.}\label{tab:path-core-certificates}
		\begin{tabular}{@{}lll@{}}
			\toprule
			Construction & Deleted edges $X$ & Three conflicting pairs\\
			\midrule
			(a), $q=5$ & $e_2,e_3,h$ & $[a_0,h];[e_1,e_3];[e_2,b_0]$\\
			(b), $q=6$ & $e_2,e_3,h$ & $[e_1,h];[e_3,e_5];[e_2,e_4]$\\
			(b), $q=7$ & $e_3,e_4,h$ & $[e_1,h];[e_3,e_5];[e_4,e_6]$\\
			(c), $q=7$ & $e_3,e_4,h$ & $[e_5,h];[e_1,e_3];[e_4,e_6]$\\
			(d), $q=8$ & $e_3,e_4,e_5$ & $[e_5,e_7];[e_1,e_3];[e_4,e_6]$\\
			(e), $q=6$ & $e_2,e_3,e_5$ & $[e_5,h];[e_1,e_3];[e_2,e_4]$\\
			(e), $q=7$ & $e_2,e_3,e_6$ & $[e_1,e_3];[e_6,h];[e_2,e_4]$\\
			(e), $q=8$ & $e_3,e_4,e_7$ & $[e_5,e_7];[e_4,h];[e_1,e_3]$\\
			\bottomrule
		\end{tabular}
	\end{table}
	
	\FloatBarrier
	The family $\mathcal R_4$ was verified in Corollary~\ref{cor:m3-cyclomatic}. For $\mathcal R_5$, label the five cycle edges consecutively as $e_1,\ldots,e_5$, with the attachment vertex common to $e_1,e_5$. Delete $e_2,e_3,e_4$; the two remaining edges form a star at the attachment vertex. The conflict graph of the five cycle edges contains the cycle $e_1e_3e_5e_2e_4e_1$, and Lemma~\ref{lem:leaf-extension} applies.
	
	Finally, consider $\mathcal S_4$. Write its branches as $c x_{i,1}\cdots x_{i,d_i}$, where $d_i\in\{2,3\}$ and $1\leq i\leq3$. Delete the three edges $cx_{i,1}$. The remaining edges form three disjoint stars, centred at $x_{i,d_i-1}$; all additional leaves are added at these centres. With indices read modulo three, the three pairs $[cx_{i,1},x_{i+1,1}x_{i+1,2}]$ are distinct conflicting pairs, using $cx_{i+1,1}$ as a witnessing edge. This gives the $3K_2$ certificate for all four branch-length patterns at once.
\end{proof}
\section{Parameterized consequences}\label{sec:parameterized}
Connectedness is unnecessary in this section, and an edgeless input has optimum zero. We use $t$ for a required number of selected edges and $k$ for the permitted number of exclusions, so the two thresholds are $\rhooe(G)\geq t$ and $\eta(G)=m-\rhooe(G)\leq k$.

\subsection{The target cardinality as parameter}
\label{subsec:solution-parameter}

For an integer $t\geq0$, let \textnormal{\textsc{Edge Open Packing}} ask whether $\rhooe(G)\geq t$, with parameter $t$. A characterization for each fixed $t$ does not by itself imply fixed-parameter tractability. We use the standard
$\mathrm{W}[1]$-completeness of \textnormal{\textsc{Independent Set}}, parameterized
by the requested cardinality~\cite{DowneyFellows1995W1}.

We first record the particular consequence of the product results of
Bre\v{s}ar et al.~\cite[Theorems 3.1 and 3.2]{BRESAR2025} that is needed
for the reduction. For a graph $H$ and an integer $q\geq2$, the graph
$H[K_q]$ replaces every vertex $v$ of $H$ by a clique $Q_v$ of order
$q$, and makes $Q_u$ complete to $Q_v$ precisely when $uv\in E(H)$.

\begin{proposition}\label{prop:clique-blowup-equality}
	For every graph $H$ and integer $q\geq2$,
	$\nuI(H[K_q])=\rhooe(H[K_q])=\alpha(H)$.
\end{proposition}
\begin{proof}
	The leaves of the star components induced by an EOP form an independent
	set whose cardinality is the number of selected edges; for a component
	$K_2$, designate either endpoint as its leaf. Consequently,
	$\rhooe(F)\leq\alpha(F)$ for every graph $F$. An independent set of
	$H[K_q]$ contains at most one vertex from each $Q_v$, and its occupied
	cliques correspond to an independent set of $H$. Hence
	$\alpha(H[K_q])=\alpha(H)$. Conversely, choose a maximum independent
	set $I$ of $H$ and one edge inside $Q_v$ for each $v\in I$. These edges
	form an induced matching of cardinality $\alpha(H)$. Combining this
	construction with $\nuI\leq\rhooe\leq\alpha$ proves the claim.
\end{proof}

The proposition is included with a direct proof for completeness. It is
a known consequence of the lexicographic-product results
in~\cite{BRESAR2025}; the parameterized consequence below also uses
their reduction principle.

\begin{theorem}\label{thm:solution-W1}
	\textnormal{\textsc{Edge Open Packing}}, parameterized by $t$, is
	$\mathrm{W}[1]$-complete. Hardness holds even on connected triangular
	graphs of diameter two that have two universal vertices and satisfy
	$\nuI(G)=\rhooe(G)$. Here, a graph is triangular if every edge belongs
	to a triangle.
\end{theorem}
\begin{proof}
	For membership, construct the edge-conflict graph $J(G)$. The equivalence
	$\rhooe(G)=\alpha(J(G))$ gives a parameter-preserving reduction to
	\textnormal{\textsc{Independent Set}}. An adjacency matrix of $G$ allows the
	conflict relation between any two edges to be checked in constant time,
	so $J(G)$ is constructed in $O(n^2+m^2)$ time.
	
	For hardness, take an instance $(H,t)$ of \textnormal{\textsc{Independent Set}}; we
	may restrict to nonempty graphs $H$. Add a vertex $x$ with no neighbour
	in $H$, then add a vertex $z$ adjacent to every vertex of
	$H\cup\{x\}$. Call the resulting graph $H'$, and set $G=H'[K_2]$.
	Since an independent set containing $z$ has size one, while
	$H\cup\{x\}$ has independence number $\alpha(H)+1$, we have
	$\alpha(H')=\alpha(H)+1$. Proposition~\ref{prop:clique-blowup-equality}
	therefore gives $\rhooe(G)=\nuI(G)=\alpha(H)+1$. Thus
	$\alpha(H)\geq t$ if and only if $\rhooe(G)\geq t+1$.
	The new parameter is $t+1$, and the construction is polynomial.
	
	Both vertices of $Q_z$ are universal in $G$. Every vertex of $Q_x$ is
	nonadjacent to every vertex in a clique corresponding to an original
	vertex of $H$, so $G$ has diameter exactly two. Each edge inside a
	clique $Q_v$ lies in a triangle using a vertex of a neighbouring clique:
	every vertex of $H'$ has a neighbour. An edge between two cliques lies
	in a triangle with the second vertex of either clique. Therefore $G$
	is triangular, as required.
\end{proof}

In particular, an algorithm with running time $f(t)n^{O(1)}$ for general
graphs would imply $\mathrm{FPT}=\mathrm{W}[1]$. This remains true even
when equality with induced matching is guaranteed. Equality alone does
not make either optimum easy to compute. Theorem~\ref{thm:solution-W1}
is a parameterized refinement of the triangular-graph NP-completeness
consequence in~\cite[Corollary 3.3]{BRESAR2025}, rather than a claim
that the underlying product identity is new.

\begin{proposition}\label{prop:target-XP}
	For every fixed $t\geq1$, deciding whether $\rhooe(G)\geq t$ takes
	$O(n^2+\binom{m}{t}t^2)$ time, with the usual immediate rejection when
	$t>m$. In particular, the problem belongs to XP.
\end{proposition}
\begin{proof}
	After constructing an adjacency matrix, enumerate the $t$-subsets of
	$E(G)$. A set is an EOP exactly when none of its $\binom{t}{2}$ pairs
	conflict. Each pair requires only a constant number of adjacency tests.
	An EOP of cardinality at least $t$ contains an EOP of cardinality exactly
	$t$, so this procedure is correct. Since $m\leq\binom{n}{2}$, its
	running time is $O(n^2+n^{2t}t^2)$.
\end{proof}

The same observation decides $\rhooe(G)\leq t$ by testing the absence
of an EOP of cardinality $t+1$. Testing both thresholds gives an XP
algorithm for the exact-value question $\rhooe(G)=t$. These statements clarify the algorithmic meaning of the structural characterization for arbitrary $t$ without asserting an FPT algorithm in that parameter.

\subsection{Adding maximum degree to the target parameter}
\label{subsec:degree-target}

The preceding hardness uses unbounded degrees. A direct kernelization
becomes available when the maximum degree $\Delta=\Delta(G)$ is also
included in the parameter.

\begin{theorem}\label{thm:degree-target-kernel}
	\textnormal{\textsc{Edge Open Packing}}, parameterized by $t+\Delta$, has a kernel
	with $O(t(\Delta+1)^2)$ vertices and edges. It is solvable in
	$2^{O(t(\Delta+1)^2)}+n^{O(1)}$ time. In particular, for every fixed
	maximum degree, the problem is fixed-parameter tractable in $t$.
\end{theorem}
\begin{proof}
	For $t=0$, return the empty graph with target zero. For $t\geq1$, reject immediately if $t>m$; this also handles edgeless inputs. We now show that $\Delta(J(G))\leq2\Delta^2$. Fix an edge $uv$.
	If an edge $f$ conflicts with $uv$, a common edge joins $u$ or $v$ to
	an endpoint of $f$. Hence, $f$ is incident with a vertex of
	$N(u)\cup N(v)$. The number of such edges is at most
	$\sum_{w\in N(u)\cup N(v)}d(w)\leq2\Delta^2$; counting an edge more
	than once only increases this upper bound.
	
	In a graph of maximum degree at most $D$, a maximal independent set
	has at least $|V|/(D+1)$ vertices: the closed neighbourhoods of its
	vertices cover the graph and each has size at most $D+1$. Applying
	this fact to $J(G)$, with $D=2\Delta^2$, shows that
	$m>(t-1)(2\Delta^2+1)$ implies $\rhooe(G)\geq t$. In this case we
	return a fixed yes-instance.
	
	Otherwise, $m\leq(t-1)(2\Delta^2+1)$. Delete the isolated vertices of
	$G$, which do not affect any edge packing. The remaining graph has at
	most $2m\leq2(t-1)(2\Delta^2+1)$ vertices and the same edge set. It
	is an equivalent instance of the same problem, and its maximum degree
	has not increased. A fixed yes- or no-instance can be used for the
	immediately decided cases, so the asserted kernel bound follows for
	all inputs. Enumerating the edge subsets of the kernel and checking
	the conflict condition proves the stated FPT running time.
\end{proof}

This kernel is a direct application of bounded-degree independent set
to the edge-conflict graph. It complements the solution-size hardness
and the separate deficit parameter considered below.
\subsection{A direct reduction to Vertex Cover}

An excluded edge of $G$ is a selected vertex of $J(G)$. For $X\subseteq E(G)$, the edges $E(G)\setminus X$ form an EOP in the original graph $G$ precisely when $X$ meets every conflict pair. Thus $\eta(G)=\tau(J(G))$, and the budget $k$ is preserved when we pass from the deficit problem to \textnormal{\textsc{Vertex Cover}}. We can use the established algorithms and reduction rules for that problem~\cite{CyganEtAl2015,BussGoldsmith1993}.

\begin{theorem}\label{thm:deficit-fpt}
	Given a graph $G$ on $n$ vertices and $m$ edges and an integer
	$k\geq0$, one can decide whether $\rhooe(G)\geq m-k$ and construct
	such an EOP, when one exists, in
	$O(2^k(n^2+m^2))$ time. In terms of $n$ alone, the running time is
	$O(2^k n^4)$.
\end{theorem}
\begin{proof}
	Construct an adjacency matrix of $G$. For each pair of edges,
	inspect the at most four pairs consisting of an endpoint of each
	edge. A conflict exists precisely when one of these pairs is a
	third edge, distinct from the two edges under consideration. Each
	test takes constant time, and $J(G)$ can be constructed
	in $O(n^2+m^2)$ time.
	
	Apply the usual two-way branching rule to $J(G)$~\cite{CyganEtAl2015}. If $xy$ is still uncovered, every feasible cover must select at least one of its endpoints. The two recursive calls select $x$ and $y$, respectively, delete the selected vertex together with its incident edges, and spend one unit of budget. An edgeless residual graph completes a cover; an uncovered edge with no budget left rejects the branch. The recursion has depth at most $k$ and at most $2^{k+1}-1$ nodes. Processing a node in $O(m^2)$ time gives the claimed bound, and removing the resulting cover from $E(G)$ yields the desired EOP. Since $m=O(n^2)$, the bound is also $O(2^k n^4)$.
\end{proof}

For a fixed $k$, the search examines at most $2^{k+1}-1$ branches,
rather than all $\binom{m}{k}$ sets of omitted edges. To recognize
the exact value $\rhooe(G)=m-k$ for $k\geq1$, run the decision
algorithm with budgets $k$ and $k-1$; the first answer must be yes
and the second no. The running time has the same asymptotic bound.

\subsection{A graph-preserving quadratic kernel for the deficit}
\label{subsec:graph-deficit-kernel}

Consider the following parameterized decision problem: given a graph $G$
and a nonnegative integer $k$, decide whether
$\rho_e^o(G)\geq |E(G)|-k$. The parameter is $k$. Unlike a reduction to
Vertex Cover, the kernel below outputs an instance of the same graph
problem.

For a vertex $v$, let $L_G(v)$ be the set of its neighbours of degree one.
We call the edges $vx$ with $x\in L_G(v)$ a pendant-edge class.

\begin{lemma}\label{lem:pendant-edge-reduction}
	Let $k\geq1$. If $|L_G(v)|>k+1$, delete all but $k+1$ vertices of
	$L_G(v)$, obtaining $G'$. Then
	$|E(G)|-\rho_e^o(G)\leq k$ if and only if
	$|E(G')|-\rho_e^o(G')\leq k$.
\end{lemma}
\begin{proof}
	Let $J(G)$ be the graph whose vertices are the edges of $G$ and whose
	adjacencies are the conflicting pairs. An EOP in $G$ is an independent set
	in $J(G)$, so $|E(G)|-\rho_e^o(G)=\tau(J(G))$. The pendant edges at $v$
	are pairwise nonadjacent in $J(G)$ and have identical open
	neighbourhoods there. Indeed, a pendant edge $vx$ conflicts with an edge
	$ab$ not incident with $v$ exactly when $v$ is adjacent to at least one
	of $a,b$; it conflicts with no edge incident with $v$.
	
	Deleting a pendant edge creates or destroys no conflict between retained
	edges: a deleted pendant edge cannot serve as a third edge joining the
	endpoints of two retained edges. Thus $J(G')$ is the induced subgraph of
	$J(G)$ obtained by deleting the corresponding vertices. A vertex cover
	of $J(G)$ restricts to a vertex cover of $J(G')$. Conversely, a vertex
	cover of $J(G')$ of size at most $k$ omits at least one of the $k+1$
	retained members of the pendant-edge class. It must therefore contain
	every vertex in their common neighbourhood. The same cover covers all
	edges incident with the deleted members in $J(G)$, and hence is also a
	vertex cover of $J(G)$.
\end{proof}

\begin{lemma}\label{lem:deficit-nonpendant-bound}
	Let $G$ be a connected graph with at least one edge that is not a star. If
	$d=|E(G)|-\rho_e^o(G)$, then $d\geq1$ and $G$ has at most $2d+1$
	vertices of degree at least two.
\end{lemma}
\begin{proof}
	Choose a maximum EOP $D$, and let $q\geq1$ be the number of nontrivial
	stars in $G[V_D]$. Let $r=|V(G)\setminus V_D|$. Each of the $d$ edges
	outside $D$ has an endpoint outside $V_D$, because $G[V_D]$ has edge set
	exactly $D$. Moreover, $r\geq1$, since otherwise the connected graph
	$G$ would itself be a star. Contracting the $q$ stars gives a connected
	multigraph on $q+r$ vertices with $d$ edges. Consequently,
	$q+r\leq d+1$.
	
	Designate one centre in each star of $G[V_D]$, choosing either endpoint
	when the star is $K_2$. A vertex of degree at least two is either one of
	these $q$ centres, one of the $r$ outside vertices, or a leaf of a star
	in $G[V_D]$ incident with an edge outside $D$. There are at most $d$
	leaves of the last kind, because each edge outside $D$ has at most one
	endpoint in $V_D$. The desired number is therefore at most
	$q+r+d\leq2d+1$.
\end{proof}

\begin{remark}\label{rem:nonpendant-sharp}
	The bound $2d+1$ in Lemma~\ref{lem:deficit-nonpendant-bound} is attained for infinitely many deficits. For $r\geq1$, take the path $P_{4r+3}$ and label its $4r+2$ edges consecutively. Two edges conflict exactly when their indices differ by two. Its conflict graph is therefore the disjoint union of two paths on $2r+1$ vertices, each having vertex-cover number $r$. Thus, the EOP deficit is $d=2r$, while the original path has $4r+1=2d+1$ vertices of degree at least two.
\end{remark}

\begin{theorem}\label{thm:graph-deficit-kernel}
	For $k\geq1$, the problem of deciding whether
	$\rho_e^o(G)\geq |E(G)|-k$, parameterized by $k$, has a kernel with
	at most $3k(k+2)$ vertices and at most $(15k^2+3k)/2$ edges. If the input
	graph is connected, these bounds improve to $(2k+1)(k+2)$ vertices and
	$(2k+1)^2$ edges. The kernel can be computed in $O(n+m)$ time and is an
	instance of the same problem. The case $k=0$ is decidable in linear time
	and admits a constant-size kernel.
\end{theorem}
\begin{proof}
	For $k=0$, the answer is yes exactly when every nontrivial component is
	a star, which can be checked in linear time. Assume $k\geq1$. Delete
	isolated vertices and components that are stars. This preserves the
	deficit, because both the number of edges and the maximum EOP number
	are additive over components and a star contributes zero deficit. If no
	component remains, return a fixed yes-instance.
	
	For each remaining component $G_i$, put
	$d_i=|E(G_i)|-\rho_e^o(G_i)$. Every $d_i$ is positive. If the input is
	a yes-instance, then $\sum_i d_i\leq k$, so there are at most $k$
	remaining components. By Lemma~\ref{lem:deficit-nonpendant-bound}, their
	total number $b$ of vertices of degree at least two satisfies
	$b\leq\sum_i(2d_i+1)\leq3k$. Hence, if $b>3k$, return a fixed
	no-instance.
	
	Apply Lemma~\ref{lem:pendant-edge-reduction} to every vertex. These
	reductions leave at least two pendant neighbours whenever a class is
	reduced. They therefore preserve the set of vertices of degree at least
	two. Every remaining leaf has one of these $b$ vertices as its neighbour,
	since the star components, including isolated edges, have been removed.
	There are now at most $b(k+1)$ leaves, yielding
	$|V(G')|\leq b(k+2)\leq3k(k+2)$. Every remaining edge either joins two
	vertices of degree at least two or is pendant, and therefore
	$|E(G')|\leq\binom b2+b(k+1)\leq(15k^2+3k)/2$. The output is $(G',k)$.
	
	For a connected input that is not a star, use the stronger rejection
	bound $b>2k+1$. The same argument gives
	$|V(G')|\leq(2k+1)(k+2)$ and
	$|E(G')|\leq\binom{2k+1}2+(2k+1)(k+1)=(2k+1)^2$.
	Computing components, degrees, pendant-edge classes, and the retained
	graph requires $O(n+m)$ time. Fixed yes- and no-instances may be used
	whenever a rule decides the answer; for example, $(K_2,0)$ is a
	yes-instance and $(P_4,0)$ is a no-instance.
\end{proof}

The kernel gives an explicit algorithm for the deficit and allows a
maximum EOP to be recovered in the original graph.

\begin{corollary}\label{cor:deficit-algorithm}
	Given a graph $G$ and $k\geq0$, one can decide whether
	$|E(G)|-\rho_e^o(G)\leq k$ in deterministic
	$O(n+m+2^k(k+1)^4)$ time. Within the same bound, one can decide whether
	$|E(G)|-\rho_e^o(G)=k$. If the deficit is at most $k$, the algorithm
	can determine its exact value and output a maximum EOP of the original
	graph. Consequently, for every fixed $k$, recognition of both
	$\rho_e^o(G)\geq m-k$ and $\rho_e^o(G)=m-k$ takes $O(n+m)$ time.
\end{corollary}
\begin{proof}
	Handle $k=0$ directly by testing whether every nontrivial component is a
	star. For $k\geq1$, compute the kernel from
	Theorem~\ref{thm:graph-deficit-kernel}. An immediate no-answer from a
	kernel rejection rule means that the deficit is greater than $k$.
	If all components were discarded as stars, their full edge set is already a maximum EOP, with deficit zero. Otherwise, let $G'$ denote the retained graph. It has $O(k^2)$ vertices
	and edges. Construct an adjacency matrix of $G'$ and then examine all
	pairs of its edges to construct $J(G')$. Each conflict is decided by at
	most four adjacency queries, so this takes $O(k^4)$ time and space.
	
	Run the Vertex Cover search from Theorem~\ref{thm:deficit-fpt} on $J(G')$ for budgets $j=0,\ldots,k$, stopping at the first successful budget. Its recursion for budget $j$ has $O(2^j)$ nodes. Each node needs at most $O(k^4)$ time to inspect the kernel conflict graph, and the total over all budgets is $O(\sum_{j=0}^k2^j k^4)=O(2^k k^4)$.
	
	The pendant reduction used for budget $k$ preserves every threshold
	$j\leq k$. Indeed, a cover of size at most $j$ omits at least one of the
	$k+1$ retained members of each reduced pendant-edge class, and therefore
	contains the common neighbours of that class. The extension argument
	of Lemma~\ref{lem:pendant-edge-reduction} works without change. Removing
	star components also preserves every threshold. Thus, the first
	successful budget is the exact original deficit, provided this value
	is at most $k$; if no budget succeeds, the original deficit exceeds
	$k$. This determines both the at-most and exact-value questions.
	
	Let $C'$ be a minimum vertex cover found in $J(G')$, where $|C'|\leq k$,
	and identify its vertices with the corresponding retained original
	edges. Reverse the pendant reductions one at a time. At each step, the
	same set $C'$ remains a vertex cover, by the extension argument in
	Lemma~\ref{lem:pendant-edge-reduction}. Thus, $C'$ covers all conflicts
	in $G$; restoring star components creates no conflicts. Hence
	$E(G)\setminus C'$ is an EOP of
	size $m-|C'|$. Since the original and retained deficits agree at this
	value, this EOP is maximum. Restoring its edges and reporting the answer
	requires $O(n+m)$ time.
\end{proof}

\begin{corollary}\label{cor:finite-deficit-cores}
	For every fixed integer $k\geq1$, there is an effectively computable finite collection $\mathcal K_k$ such that a connected graph $G$ satisfies $\rhooe(G)=m-k$ if and only if it can be obtained from some $H\in\mathcal K_k$ by adding arbitrary numbers of leaves at vertices that have exactly $k+1$ pendant neighbours in $H$. Zero additions are allowed. Every $H\in\mathcal K_k$ has at most $(2k+1)(k+2)$ vertices and at most $(2k+1)^2$ edges.
\end{corollary}
\begin{proof}
	Let $\mathcal K_k$ contain all connected graphs $H$ within the stated size bounds for which $\eta(H)=k$ and every vertex has at most $k+1$ pendant neighbours. This collection is finite and computable: enumerate the bounded graphs and compute the vertex-cover number of each conflict graph.
	
	Given connected $G$ with $\eta(G)=k$, cap every pendant-edge class at $k+1$. The resulting graph $H$ remains connected and, by the proof of Corollary~\ref{cor:deficit-algorithm}, preserves every deficit threshold at most $k$. In particular, it has deficit at most $k$ and not at most $k-1$, so $\eta(H)=k$. The connected kernel bounds place $H$ in $\mathcal K_k$. All deleted leaves belonged to classes that now have exactly $k+1$ members, giving the asserted construction.
	
	Conversely, add leaves only at those vertices of $H$ that already have $k+1$ pendant neighbours. Capping the enlarged classes restores $H$, and the pendant reduction preserves both thresholds $k$ and $k-1$. Hence, the expanded graph has deficit exactly $k$. As $k+1\geq2$, eligible vertices are nonleaves, so additions at one eligible vertex do not change which vertices belong to the other pendant classes.
\end{proof}

The collection $\mathcal K_k$ is computable for each fixed $k$, although enumerating it may be costly. For deficit three, Section~\ref{EOP_sec_3} gives the shorter construction directly.

\subsection{Treewidth as a structural parameter}
\label{subsec:treewidth-parameter}
Das and Santra~\cite{DasSantra2026ATFree} already noted the treewidth consequence for EOP. We include the statement to specify the required input and the objective being optimized. It follows from the general theorem of Courcelle and Mosbah~\cite{CourcelleMosbah1993}.

\begin{proposition}\label{prop:treewidth-FPT}
	Given a tree decomposition of $G$ of width $w$ with $O(n)$ bags, a maximum EOP can be computed in $f(w)(n+m)$ time for a computable function $f$.
\end{proposition}
\begin{proof}
	In MSO$_2$, edge variables and vertex variables are linked by incidence. For edges $e\neq f$, the conflict relation is defined by an edge $g\notin\{e,f\}$ and vertices $u\neq v$ such that $u$ is incident with both $e$ and $g$, and $v$ with both $f$ and $g$. We can therefore describe a feasible edge set $D$ by forbidding this relation for every pair of its members. The formula has constant length, independent of $G$. Optimizing $|D|$ under that formula is covered by the monadic second-order optimization theorem~\cite{CourcelleMosbah1993}. Applied to the supplied decomposition, it returns both the optimum value and a corresponding edge set in $f(w)(n+m)$ time.
\end{proof}

Table~\ref{tab:parameter-results} summarizes the parameterized results. The direct reduction to Vertex Cover also permits sharper algorithms such as those of Chen et al.~\cite{ChenKanjXia2010}; elementary branching suffices for our stated bounds.

\begin{table}[htbp!]
	\centering
	\small
	\renewcommand{\arraystretch}{1.25}
	\caption{Parameterized consequences for edge open packing. Here $k$ is the supplied upper bound on $m-\rhooe(G)$ in the deficit decision problem; $w$ is the width of a supplied tree decomposition. All stated kernels output graphs for the EOP problem itself.}
	\label{tab:parameter-results}
	\begin{tabular}{@{}>{\raggedright\arraybackslash}p{.19\linewidth}>{\raggedright\arraybackslash}p{.72\linewidth}@{}}
		\toprule
		Parameter & Result\\
		\midrule
		Target $t$ & $\mathrm{W}[1]$-complete; XP time $O(n^2+n^{2t}t^2)$ for the at-least question.\\
		$t+\Delta(G)$ & $O(t(\Delta+1)^2)$-vertex and edge kernel; fixed-parameter tractable.\\
		Deficit bound $k$ & $3k(k+2)$-vertex kernel for general graphs, and $(2k+1)(k+2)$ for connected graphs; $O(k^2)$ edges.\\
		Deficit bound $k$ & Deterministic time $O(n+m+2^k(k+1)^4)$ for at-most or exact deficit; maximum EOP reconstruction when the deficit is at most $k$.\\
		Fixed deficit $k$ & $O(n+m)$ recognition and reconstruction.\\
		Treewidth $w$ & Time $f(w)(n+m)$, by MSO$_2$ optimization.\\
		\bottomrule
	\end{tabular}
\end{table}
\section{Graphs with equal induced matching and edge open packing numbers}
\label{sec:equal-parameters}

Write $\nu_I(G)$ for the induced matching number of $G$. Every induced
matching is an edge open packing, so $\nu_I(G)\leq\rho_e^o(G)$.
Chelladurai et al.~\cite{chelladurai2022edge} asked for a characterization
of the graphs attaining equality. Bre\v{s}ar et al.~\cite{BRESAR2025}
settled this question for trees. Their construction starts with disjoint
spiders, joins them by edges incident with spider centres, and retains at
least two leaves of every spider. The general question appears in Section~4, item~(2)(i), of~\cite{chelladurai2022edge}, using $P_2^\prime(G)$ for the induced matching number. The next result removes the requirement
that the resulting graph be a tree.

A \emph{spider with $k$ arms} is obtained from $K_{1,k}$ by subdividing each
edge exactly once. We denote its centre by $c$, its support vertices by
$a_1,\ldots,a_k$, and its leaves by $b_1,\ldots,b_k$; its edges are
$ca_j$ and $a_jb_j$, for $1\leq j\leq k$.

\begin{theorem}\label{thm:spider-extension-equality}
	Let $G$ be obtained from vertex-disjoint spiders
	$S_1,\ldots,S_h$, where $S_i$ has centre $c_i$ and $k_i\geq2$ arms, by
	adding edges subject to the following two conditions:
	\begin{enumerate}[label=\textup{(\roman*)}]
		\item every added edge has at least one endpoint in
		$C=\{c_1,\ldots,c_h\}$;
		\item at least two original leaves of each $S_i$ have degree one in $G$.
	\end{enumerate}
	Then $\nu_I(G)=\rho_e^o(G)=\sum_{i=1}^h k_i$. Moreover, the original
	pendant edges of the spiders form the unique maximum induced matching
	of $G$.
\end{theorem}
\begin{proof}
	Let $M$ consist of the original pendant edges, and put
	$K=|M|=\sum_i k_i$. Since every added edge meets $C$, the graph
	$G-C$ consists exactly of the edges of $M$. Thus $M$ is an induced
	matching and $K\leq\nu_I(G)\leq\rho_e^o(G)$.
	
	Let $D$ be any edge open packing, and let $W$ be its set of endpoints.
	The graph $G[W]$ is a disjoint union of stars and has edge set $D$.
	Write $q$ for its number of components and $c=|C\cap W|$. For
	$j\in\{0,1,2\}$, let $p_j$ count the edges of $M$ with exactly $j$
	endpoints in $W$. Counting the vertices and edges of this star forest
	gives $|D|=2p_2+p_1+c-q$ and $K=p_0+p_1+p_2$.
	
	Choose a centre of each component of $G[W]$, arbitrarily for a component
	isomorphic to $K_2$. A component contains at most one edge of $M$, since
	$M$ is a matching, and any two edges of a star meet. If a component
	contains an edge of $M$, its chosen centre lies outside $C$. Consequently,
	at most $q-p_2$ vertices of $C\cap W$ are chosen star centres.
	
	Now let $c_i\in C\cap W$ be a vertex not chosen as a star centre. It has
	exactly one neighbour in $W$. Choose two original leaves $b_{i,r}$ and
	$b_{i,s}$ that retain degree one in $G$, with corresponding support
	vertices $a_{i,r}$ and $a_{i,s}$. At most one of these support vertices
	belongs to $W$. If $a_{i,j}\notin W$, then $b_{i,j}\notin W$ as well:
	the only neighbour of $b_{i,j}$ is $a_{i,j}$, whereas every vertex of
	$W$ is incident with an edge of $D$. Thus, at least one of these two
	edges of $M$ has no endpoint in $W$. Assign one such edge to $c_i$.
	These assignments are distinct, since the spiders are vertex-disjoint.
	Hence at most $p_0$ vertices of $C\cap W$ are star leaves. We obtain
	$c\leq q-p_2+p_0$, and therefore
	$|D|=2p_2+p_1+c-q\leq p_0+p_1+p_2=K$.
	This proves the equality of the two parameters.
	
	For uniqueness, let $N$ be an induced matching and use the same notation
	with $W$ its endpoint set. Every vertex of $C\cap W$ now has exactly
	one neighbour in $W$, so the preceding assignment gives $p_0\geq c$.
	Here $q=|N|$, and $r=q-p_2$ is the number of edges of $N$ outside $M$.
	The counting identity above yields
	$|N|=K-p_0+c-r\leq K-r$. Thus, any induced matching using an edge
	outside $M$ has size less than $K$, and $M$ is the unique maximum one.
\end{proof}

Theorem~\ref{thm:spider-extension-equality} includes every nontrivial tree
in the characterization of Bre\v{s}ar et al.~\cite{BRESAR2025} except
$P_2$, which satisfies the equality separately. Their two conditions
on added edges are retained, but neither the number of added edges
nor acyclicity is required here. In particular, edges may be added
between centres, from a centre to a support vertex of another spider,
or from a centre to a leaf, provided that two original leaves of each
spider retain degree one. The graph may be disconnected.

\begin{corollary}\label{cor:arbitrary-core-equality}
	Let $H$ be any graph. For each $v\in V(H)$, choose an integer $k_v\geq2$
	and add new vertices $a_{v,j},b_{v,j}$ and edges
	$va_{v,j},a_{v,j}b_{v,j}$ for $1\leq j\leq k_v$, with all new vertices
	distinct. The resulting graph $G$ satisfies
	$\nu_I(G)=\rho_e^o(G)=\sum_{v\in V(H)}k_v$.
\end{corollary}
\begin{proof}
	Apply Theorem~\ref{thm:spider-extension-equality}, taking the edges of
	$H$ as the added edges between spider centres.
\end{proof}

If all $k_v=k$, this is the rooted product $H\circ_c S_k$, where the
root $c$ is the centre of the spider $S_k$, and the common value is
$k|V(H)|$. Thus, the equality holds for a class containing every graph
as an induced subgraph. This gives a particular exact formula for the
rooted-product problem posed by Bre\v{s}ar et al.~\cite{BRESAR2025};
it does not settle their question for arbitrary rooted graphs.

\begin{figure}[htbp]
	\centering
	\begin{tikzpicture}[eop graph]
		\eopv{c1}{-1.6}{0}\eopv{c2}{1.6}{0}\eopv{c3}{0}{2.2}
		\draw[dashed] (c1)--(c2)--(c3)--(c1);
		\foreach \c/\x/\y/\dx/\dy/\id in {c1/-1.6/0/-1/-.6/1,c1/-1.6/0/-.9/.65/2,c2/1.6/0/1/-.6/3,c2/1.6/0/.9/.65/4,c3/0/2.2/-.65/.85/5,c3/0/2.2/.65/.85/6}{
			\node[eop vertex] (a\id) at (\x+\dx,\y+\dy) {};
			\node[eop vertex] (b\id) at (\x+1.8*\dx,\y+1.8*\dy) {};
			\draw (\c)--(a\id);\draw[line width=1.1pt](a\id)--(b\id);}
		\node[eop label,below=3pt] at(c1){$c_1$};
		\node[eop label,below=3pt] at(c2){$c_2$};
		\node[eop label,below=3pt] at(c3){$c_3$};
	\end{tikzpicture}
	\caption{Three spiders with two arms each, joined through their centres by the dashed triangle edges. The six thick edges form the unique maximum induced matching, and the EOP number is also six. The centre graph may be arbitrary.}
	\label{fig:equality-spiders}
\end{figure}

When the spider decomposition is supplied, the common optimum is constructed in $O(n+m)$ time: output the original pendant edges and verify the two conditions by scanning the edges and leaf degrees. This does not assert a recognition algorithm for all equality graphs.

The number of retained leaves is essential even when every original
spider has two arms. Take two such spiders, with centres $c,d$ and arms
$ca_1b_1,ca_2b_2$ and $dx_1y_1,dx_2y_2$, and add the edge $db_1$.
Only one original leaf of the first spider retains degree one. The four
original pendant edges form an induced matching, and an induced matching
of size five would cover all ten vertices and force the whole connected
graph to be a matching. Thus, its induced matching number is four.
However, $\{ca_2,a_2b_2,db_1,dx_1,dx_2\}$ is an edge open packing
of size five. A spider with one arm is $P_3$, whose induced matching and edge open packing
numbers are one and two, respectively.

\begin{proposition}\label{prop:hereditary-equality}
	For a graph $G$, the equality $\nu_I(F)=\rho_e^o(F)$ holds for every
	induced subgraph $F$ of $G$ if and only if $G$ is a disjoint union of
	cliques.
\end{proposition}
\begin{proof}
	A graph is a disjoint union of cliques if and only if it has no induced
	$P_3$. Since the two parameters differ on $P_3$, the stated equality for
	all induced subgraphs forces this structure. Conversely, in a disjoint
	union of cliques both parameters equal the number of components
	containing an edge. This also holds for each induced subgraph.
\end{proof}

\section{Conclusion}\label{EOP_sec_4}
Prescribing the EOP number and bounding its deficit lead to two useful structural descriptions. Our criterion handles every positive optimum, while deficit three is captured by a star-extension rule and four conflict certificates. The explicit constructions and their certificates are verified in the main text. In the connected case, deficit three forces a vertex cover of size at most four and at most two independent cycles. The graphs with two independent cycles have the theta-graph description given above.

The target size and the deficit have different parameterized behavior. The target-size problem is $\mathrm{W}[1]$-complete, whereas the deficit admits a linear-time quadratic graph kernel and a deterministic $O(n+m+2^k(k+1)^4)$-time algorithm. Consequently, every fixed-deficit recognition problem, including $\rhooe(G)=m-3$, is solvable in linear time. The combined parameter $t+\Delta(G)$ also admits a kernel, while treewidth gives fixed-parameter tractability through standard logical optimization.

The equality theorem permits arbitrary centre-incident extensions of disjoint spiders as long as two original leaves per spider remain pendant. It extends the equality construction for trees to a class containing arbitrary induced graph cores. A characterization of all graphs with $\nuI(G)=\rhooe(G)$ remains open. Further questions include whether the quadratic deficit kernel can be improved to a linear-size graph kernel, which other rooted graphs admit an exact EOP formula for every core, and how to obtain small explicit template lists for larger fixed deficits beyond the computable finite-core construction.

\bibliographystyle{plain}
\bibliography{EOP_Bib}

\end{document}